\documentclass[11pt]{article}
\usepackage{url}
\begin{document}
%
\font\tenrm=cmr10
\font\ninerm=cmr9
\font\eightrm=cmr8
\font\sevenrm=cmr7
\font\ss=cmss10 
\font\smallcaps=cmcsc10 
%
\overfullrule=0pt
%
%
%
%
\newcommand{\mem}[2]{\chi_{#1}^{#2}}
\newcommand{\implies}{\rightarrow}
\newcommand{\el}{\ell }
\newcommand{\eltwo}{{\small \mathcal{L} }}
\newcommand{\sigp}[1]{\Sigma_{#1}^{{\rm p}}}
\newcommand{\pip}[1]{\Pi_{#1}^{{\rm p}}}
\newcommand{\sigsp}[1]{\Sigma_{#1}^{{\rm p, SPARSE}}}
\newcommand{\pisp}[1]{\Pi_{#1}^{{\rm p, SPARSE}}}
\newcommand{\Reone}{R_e^1}
\newcommand{\Retwo}{R_e^2}
\newcommand{\Rethree}{R_e^3}
\newcommand{\Reonesub}[1]{R_{\ang{e,{#1}}}^1}
\newcommand{\Retwosub}[1]{R_{\ang{e,{#1}}}^2}
\newcommand{\dom}{{\rm dom }}
\newcommand{\ran}{{\rm ran }}
\newcommand{\nd}{{\rm nd }}
\newcommand{\rd}{{\rm rd }}
\newcommand{\ith}{{i}^{\rm th}}
\newcommand{\Prob}{{\rm Pr } }
\renewcommand{\Pr}{{\rm Pr } }
\newcommand{\pr}{{\rm Pr } }
\newcommand{\DTIME}{{\rm DTIME } }
\newcommand{\NTIME}{{\rm NTIME } }
\newcommand{\DSPACE}{{\rm DSPACE } }
\newcommand{\NSPACE}{{\rm NSPACE } }
%
%

\newcommand{\PH}{{\rm PH}}
\newcommand{\FML}{{\rm FML}}
\newcommand{\ED}{{\rm ED}}
\newcommand{\Max}{{\rm MAX}}
\newcommand{\Med}{{\rm MED}}
\newcommand{\Sort}{{\rm SORT}}
\newcommand{\Sel}{{\rm SEL}}

\newcommand{\mink}{{\rm MIN}_k}
\newcommand{\maxn}{{\rm MAX}_n}
\newcommand{\maxk}{{\rm MAX}_k}
\newcommand{\maxnm}{{\rm MAX}_{n-1}}
\newcommand{\selin}{{\rm SEL}_n^i}
\newcommand{\selina}{{\rm SELALL}_n^i}
\newcommand{\secn}{{\rm SEC}_n}
\newcommand{\medn}{{\rm MED}_n}
\newcommand{\sortn}{{\rm SORT}_n}
\newcommand{\sortk}{{\rm SORT}_k}
\newcommand{\hilon}{{\rm HI\hbox{-}LO}_n}

\newcommand{\Abar}{\overline A}
\newcommand{\Bbar}{\overline B}
\newcommand{\COFbar}{{\overline {COF}}}
\newcommand{\Kbar}{\overline K}
\newcommand{\Kleenestar}{{\textstyle *}}
\newcommand{\Kstar}{{\textstyle *}}
\newcommand{\NE}{ {\rm NE} }
\newcommand{\SPARSE}{ {\rm SPARSE} }
\newcommand{\TOTbar}{{\overline {TOT}}}
\newcommand{\Uei}{U_{\ang{e,i}}}
\newcommand{\ang}[1]{\langle#1\rangle}
\newcommand{\co}{ {\rm co} }
\newcommand{\cvg}{\downarrow}
\newcommand{\dvg}{\uparrow}
\newcommand{\ei}{\ang{e,i}}
\newcommand{\fuinv}[2]{{#1}^{\hbox{-}1}}
\renewcommand{\gets}{\leftarrow}
\newcommand{\goes}{\rightarrow}
\newcommand{\jo}{\oplus}
\newcommand{\kstar}{{\textstyle *}}
\newcommand{\move}{\vdash}
\newcommand{\penpns}{\hbox{${\rm P}={\rm NP}$}}
\newcommand{\penpq}{\hbox{P$=$?NP\ }}
\newcommand{\penp}{\hbox{P$=$NP\ }}
\newcommand{\phivec}[1]{{\ang {\phi_1,\ldots,\phi_{#1}}}}
\newcommand{\ppolyA}{{{\rm P}^A/{\rm poly}}}
\newcommand{\nppolyA}{{{\rm NP}^A/{\rm poly}}}
\newcommand{\xibar}{{\bar x_i}}
\newcommand{\xjbar}{{\bar x_j}}
\newcommand{\xkbar}{{\bar x_k}}
\newcommand{\yei}{y_{\ang{e,i}}}
\newcommand{\set}[2]{ \{ {#1}_1,\ldots,{#1}_{#2} \} }
\newcommand{\m}{{\rm m}}
\newcommand{\ctt}{{\rm ctt}}
\newcommand{\wtt}{{\rm wtt}}
\newcommand{\T}{{\rm T}}
\newcommand{\btt}{{\rm btt}}
\newcommand{\bwtt}{{\rm bwtt}}
\newcommand{\ktt}{{k\rm\mbox{-}tt}}
\newcommand{\ntt}{{n\rm\mbox{-}tt}}
\newcommand{\kwtt}{{k\rm\mbox{-}wtt}}
\newcommand{\nwtt}{{n\rm\mbox{-}wtt}}
\newcommand{\kT}{{k\rm\mbox{-}T}}
\newcommand{\onett}{{1\rm\mbox{-}tt}}
\newcommand{\leonett}{{\le_{\onett}}}
\newcommand{\eqonett}{{\equiv_{\onett}}}
\newcommand{\pTE}{\equiv_{\rm T}^p}
\newcommand{\UNIQTSP}{\rm UNIQTSP}
\newcommand{\TAUT}{{\rm TAUT}}
\newcommand{\SAT}{{\rm SAT}}
\newcommand{\usat}{{\rm UNIQSAT}}
\newcommand{\paritysat}{{\rm PARITYSAT}}
\newcommand{\LSAT}{{\rm LSAT}}
\newcommand{\R}{{\rm R}}
\renewcommand{\P}{{\rm P}}
\newcommand{\PCP}{{\rm PCP}}
\newcommand{\NP}{{\rm NP}}
\newcommand{\BPP}{{\rm BPP}}
\newcommand{\AM}{{\rm AM}}
\newcommand{\MA}{{\rm MA}}
\newcommand{\IP}{{\rm IP}}
\newcommand{\coNP}{\rm co\hbox{-}NP}
\newcommand{\QBF}{\rm QBF}
\newcommand{\PSPACE}{\rm PSPACE}
\newcommand{\alephnot}{\aleph_0}
\newcommand{\twoalephnot}{2^{\aleph_0}}
\newcommand{\lamx}[1]{\lambda x [{#1}] }
\newcommand{\bits}[1]{\{0,1\}^{{#1}}}
\newcommand{\bit}{\{0,1\}}
\newcommand{\PF}{{\rm PF}}
\newcommand{\FP}{{\rm PF}}
\newcommand{\poly}{{\rm poly}}
\newcommand{\ppoly}{\P/\poly}
\newcommand{\nppoly}{\NP/\poly}
\newcommand{\conppoly}{\coNP/\poly}
\newcommand{\Z}{{\sf Z}}
\newcommand{\nat}{{\sf N}}
\newcommand{\rat}{{\sf Q}}
\newcommand{\real}{{\sf R}}
\newcommand{\complex}{{\sf C}}
\newcommand{\Rpos}{{\sf R}^+}
\newcommand{\F}[2]{{\rm F}_{#1}^{#2}}
\newcommand{\FK}[1]{{\rm F}_{#1}^{K}}
\newcommand{\FA}[1]{{\rm F}_{#1}^{A}}
\newcommand{\FAp}[1]{{\rm F}_{#1}^{A'}}
\newcommand{\FB}[1]{{\rm F}_{#1}^{A}}
\newcommand{\FC}[1]{{\rm F}_{#1}^{A}}
\newcommand{\FX}[1]{{\rm F}_{#1}^{X}}
\newcommand{\FY}[1]{{\rm F}_{#1}^{Y}}
\newcommand{\V}[2]{{\rm V}_{#1}^{#2}}
\newcommand{\VK}[1]{{\rm V}_{#1}^{K}}
\newcommand{\VA}[1]{{\rm V}_{#1}^{A}}
\newcommand{\VB}[1]{{\rm V}_{#1}^{B}}
\newcommand{\VC}[1]{{\rm V}_{#1}^{C}}
\newcommand{\VX}[1]{{\rm V}_{#1}^{X}}
\newcommand{\VY}[1]{{\rm V}_{#1}^{Y}}
\newcommand{\PFq}[2]{\PF_{{#1\rm\mbox{-}T}}^{#2}}
\newcommand{\Q}[2]{\P_{{#1\rm\mbox{-}T}}^{#2}}
\newcommand{\EN}[1]{{\rm EN}({#1})}
\newcommand{\ENA}[1]{{\rm EN}^A({#1})}
\newcommand{\ENX}[1]{{\rm EN}^X({#1})}
\newcommand{\FQp}[2]{\PF_{{#1\rm\mbox{-}tt}}^{#2}}
\newcommand{\Qp}[2]{\P_{{#1\rm\mbox{-}tt}}^{#2}}
\newcommand{\PARITY}[2]{{\rm PARITY}_{#1}^{#2}}
\newcommand{\parity}[2]{{\rm PARITY}_{#1}^{#2}}
\newcommand{\odd}[2]{{\rm ODD}_{#1}^{#2}}
\newcommand{\MODm}[2]{  { {\rm MOD}m}_{#1}^{#2}}
\newcommand{\SMODm}[2]{  { {\rm SMOD}m}_{#1}^{#2}}
\newcommand{\PARITYP}{\PARITY\P}
\newcommand{\MOD}{{\rm MOD}}
\newcommand{\SMOD}{{\rm SMOD}}
\newcommand{\POW}{{\rm POW}}
\newcommand{\FAC}{{\rm FAC}}
\newcommand{\POLY}{{\rm POLY}}
\newcommand{\card}[2]{\#_{#1}^{#2}}
\newcommand{\pleq}[1]{\leq_{#1}^{\rm p}}
\newcommand{\plem}{\le_{\rm m}^{\rm p}}
\newcommand{\pleT}{\le_{\rm T}^{\rm p}}
%
%
\newcommand{\xvec}[1]{\ifcase 3{#1} {\ang {x_1,x_2,x_3} } \else 
\ifcase 4{#1} {\ang{x_1,x_2,x_3,x_4}} \else {\ang {x_1,\ldots,x_{#1}}}\fi\fi}
\newcommand{\yvec}[1]{\ifcase 3{#1} {\ang {y_1,y_2,y_3} } \else 
\ifcase 4{#1} {\ang{y_1,y_2,y_3,y_4}} \else {\ang {y_1,\ldots,y_{#1}}}\fi\fi}
\newcommand{\zvec}[1]{\ifcase 3{#1} {\ang {z_1,z_2,z_3} } \else 
\ifcase 4{#1} {\ang{z_1,z_2,z_3,z_4}} \else {\ang {z_1,\ldots,z_{#1}}}\fi\fi}
\newcommand{\vecc}[2]{\ifcase 3{#2} {\ang { {#1}_1,{#1}_2,{#1}_3 } } \else
\ifcase 4{#1} {\ang { {#1}_1,{#1}_2,{#1}_3,{#1}_{4} } }
\else {\ang { {#1}_1,\ldots,{#1}_{#2}}}\fi\fi}
\newcommand{\veccd}[3]{\ifcase 3{#2} {\ang { {#1}_{{#3}1},{#1}_{{#3}2},{#1}_{{#3}3} } } \else
\ifcase 4{#1} {\ang { {#1}_{{#3}1},{#1}_{{#3}2},{#1}_{#3}3},{#1}_{{#3}4} }
\else {\ang { {#1}_{{#3}1},\ldots,{#1}_{{#3}{#2}}}}\fi\fi}
\newcommand{\xyvec}[1]{\ang{ \ang{x_1,y_1},\ang{x_2,y_2}\ldots,\ang{x_{#1},y_{#1}}}}
%
%
%
\newcommand{\veccz}[2]{\ifcase 3{#2} {\ang { {#1}_0,{#1}_2,{#1}_3 } } \else
\ifcase 4{#1} {\ang { {#1}_0,{#1}_2,{#1}_3,{#1}_{4} } }
\else {\ang { {#1}_0,\ldots,{#1}_{#2}}}\fi\fi}
\newcommand{\phive}[1]{{\phi_1,\ldots,\phi_{#1}}}
%
%
\newcommand{\xve}[1]{\ifcase 3{#1} {x_1,x_2,x_3} \else 
\ifcase 4{#1} {x_1,x_2,x_3,x_4} \else {x_1,\ldots,x_{#1}}\fi\fi}
\newcommand{\yve}[1]{\ifcase 3{#1} {y_1,y_2,y_3} \else 
\ifcase 4{#1} {y_1,y_2,y_3,y_4} \else {y_1,\ldots,y_{#1}}\fi\fi}
\newcommand{\zve}[1]{\ifcase 3{#1} {z_1,z_2,z_3} \else 
\ifcase 4{#1} {z_1,z_2,z_3,z_4} \else {z_1,\ldots,z_{#1}}\fi\fi}
\newcommand{\ve}[2]{\ifcase 3#2 {{#1}_1,{#1}_2,{#1}_3} \else
\ifcase 4#2 {{#1}_1,{#1}_2,{#1}_3,{#1}_{4}}
\else {{#1}_1,\ldots,{#1}_{#2}}\fi\fi}
\newcommand{\ved}[3]{\ifcase 3#2 {{#1}_{{#3}1},{#1}_{{#3}2},{#1}_{{#3}3}} \else
\ifcase 4#2 {{#1}_{{#3}1},{#1}_{{#3}2},{#1}_{{#3}3},{#1}_{{#3}4}}
\else {{#1}_{{#3}1},\ldots,{#1}_{{#3}{#2}}}\fi\fi}
\newcommand{\fuve}[3]{
\ifcase 3#2
{{#3}({#1}_1),{#3}({#1}_2,{#3}({#1}_3)} \else
\ifcase 4#2
{{#3}({#1}_1),{#3}({#1}_2),{#3}({#1}_3),{#3}({#1}_4)}
\else
{{#3}({#1}_1),\ldots,{#3}({#1}_{#2})}\fi\fi}
\newcommand{\fuvec}[3]{\ang{\fuve{#1}{#2}{#3}}}
\newcommand{\xse}[1]{\xve{#1}}
\newcommand{\yse}[1]{\yve{#1}}
\newcommand{\zse}[1]{\zve{#1}}
\newcommand{\fuse}[3]{\fuve{#1}{#2}{#3}}
%
%
\newcommand{\xset}[1]{\{\xve{#1}\}}
\newcommand{\yset}[1]{\{\yve{#1}\}}
\newcommand{\zset}[1]{\{\zve{#1}\}}
\newcommand{\setd}[3]{\{\ved{#1}{#2}{#3}\}}
\newcommand{\fuset}[3]{\{\fuve{#1}{#2}{#3}\}}
%
%

\newcommand{\setmathchar}[1]{\ifmmode#1\else$#1$\fi}
\newcommand{\vlist}[2]{%
	\setmathchar{%
		\compound#2\one{#2}\two
		\ifcompound
			({#1}_1,\ldots,{#1}_{#2})
		\else
			\ifcat N#2
				({#1}_1,\ldots,{#1}_{#2})
			\else
				\ifcase#2
					({#1}_0)\or
					({#1}_1)\or
					({#1}_1,{#1}_2)\or 
					({#1}_1,{#1}_2,{#1}_3)\or
					({#1}_1,{#1}_2,{#1}_3,{#1}_4)\else 
					({#1}_1,\ldots,{#1}_{#2})
				\fi
			\fi
		\fi}}
\newcommand{\xtu}[1]{\vlist{x}{#1}}
\newcommand{\ytu}[1]{\vlist{y}{#1}}
\newcommand{\ztu}[1]{\vlist{z}{#1}}
\newcommand{\btu}[1]{\vlist{b}{#1}}
\newcommand{\ptu}[1]{\vlist{p}{#1}}
\newcommand{\qtu}[1]{\vlist{q}{#1}}
\newcommand{\tup}[2]{\vlist{#2}{#1}}

\newif\ifcompound
\def\compound#1\one#2\two{%
	\def\one{#1}
	\def\two{#2}
	\if\one\two
		\compoundfalse
	\else
		\compoundtrue
	\fi}

\newcommand{\tu}[2]{(\ve{#1}{#2})}
\newcommand{\tud}[3]{(\ve{#1}{#2}{#3})}

\newcommand{\futu}[3]{(\fuve#1#2#3)}
%
%
\newcommand{\xwe}[1]{\ifcase 3{#1} {x_1\wedge x_2\wedge x_3} \else 
\ifcase 4{#1} {x_1\wedge x_2\wedge x_3\wedge x_4} \else {x_1\wedge \cdots \wedge
x_{#1}}\fi\fi}
\newcommand{\we}[2]{\ifcase 3#2 {\ang { {#1}_1\wedge {#1}_2\wedge {#1}_3 } } \else
\ifcase 4{#1} {\ang { {#1}_1\wedge {#1}_2\wedge {#1}_3\wedge {#1}_{4} } }
\else {\ang { {#1}_1\wedge \cdots\wedge {#1}_{#2}}}\fi\fi}
\newcommand{\phione}{\es'}
\newcommand{\phidoub}{\es''}
\newcommand{\phitrip}{\es'''}
\newcommand{\phiomega}{\es^{\omega}}
\newcommand{\dphione}{{\bf 0'}}
\newcommand{\dphidoub}{{\bf 0''}}
\newcommand{\dphitrip}{{\bf 0'''}}
\newcommand{\st}{\mathrel{:}}
\newcommand{\e}{\{e\}}
\newcommand{\eee}[2]{\{e\}_{#2}^{#1}}
\newcommand{\ess}{\{e\}_s}
\newcommand{\esub}[1]{\{e_{#1}\}}
\newcommand{\esubsub}[2]{{\{e_{#1}\}}_{#2}}
\newcommand{\et}{\{e\}_t}
\newcommand{\join}{\oplus}
\newcommand{\poneE}{\equiv_{\rm p}^1}
\newcommand{\mE}{\equiv_{\rm m}}
\newcommand{\pmE}{\equiv_{\rm m}^{\rm p}}
\newcommand{\phiawe}{\Phi_e^A}
\newcommand{\hpmE}{\equiv_{\rm m}^{\rm h}}
\newcommand{\iE}{\equiv_{\rm i}}
\newcommand{\TE}{\equiv_{\rm T}}
\newcommand{\hpTE}{\equiv_{\rm T}^{\rm h}}
\newcommand{\ttE}{\equiv_{{\rm tt}}}
\newcommand{\TJ}{T_j}
\newcommand{\mLE}{\le_{\rm m}}
\newcommand{\dLE}{\le_{\rm d}}
\newcommand{\cLE}{\le_{\rm c}}
\newcommand{\pLE}{\le_{\rm p}}
\newcommand{\rLE}{\le_{\rm r}}
\newcommand{\dE}{\equiv_{\rm d}}
\newcommand{\cE}{\equiv_{\rm c}}
\newcommand{\pE}{\equiv_{\rm p}}
\newcommand{\rE}{\equiv_{\rm r}}
\newcommand{\oneLE}{\le_1}
\newcommand{\oneE}{\equiv_1}
\newcommand{\oneL}{<_1}
\newcommand{\pmLE}{\le_{\rm m}^{\rm p}}
\newcommand{\hpmLE}{\le_{\rm m}^{\rm h}}
\newcommand{\hpmL}{<_{\rm m}^{\rm h}}
\newcommand{\hpmtoE}{\equiv_{\rm m}^{{\rm h}\hbox{-}{\rm to}}}
\newcommand{\hpmtoL}{<_{\rm m}^{{\rm h}\hbox{-}{\rm to}}}
\newcommand{\hpmtoLE}{\le_{\rm m}^{{\rm h}\hbox{-}{\rm to}}}
\newcommand{\hpmto}{\le_{\rm m}^{{\rm h}\hbox{-}{\rm to}}}
\newcommand{\Ba}{\hbox{{\bf a} }}
\newcommand{\Bb}{\hbox{{\bf b} }}
\newcommand{\Bc}{\hbox{{\bf c} }}
\newcommand{\Bd}{\hbox{{\bf d} }}
\newcommand{\nre}{\hbox{{\it n}-r.e.}}
\newcommand{\Bnre}{\hbox{{\bf n-r.e.}}}
\newcommand{\into}{\rightarrow}
\renewcommand{\AE}{\forall^\infty}
\newcommand{\IO}{\exists^\infty}
\newcommand{\ep}{\epsilon}
\newcommand{\es}{\emptyset}
\newcommand{\isom}{\simeq}
\newcommand{\nisom}{\not\simeq}

\newcommand{\lf}{\left\lfloor}
\newcommand{\rf}{\right\rfloor}
\newcommand{\lc}{\left\lceil}
\newcommand{\rc}{\right\rceil}
\newcommand{\Ceil}[1]{\left\lceil {#1}\right\rceil}
\newcommand{\ceil}[1]{\left\lceil {#1}\right\rceil}
\newcommand{\floor}[1]{\left\lfloor{#1}\right\rfloor}

\newcommand{\nth}{n^{th}}
\newcommand{\lecc}{\le_{\rm cc}}
\newcommand{\TLE}{\le_{\rm T}}
\newcommand{\ttLE}{\le_{\rm tt}}
\newcommand{\nttLE}{\not\le_{\rm tt}}
\newcommand{\bttLE}{\le_{\btt}}
\newcommand{\bttE}{\equiv_{\btt}}
\newcommand{\wttLE}{\le_{\wtt}}
\newcommand{\wttE}{\equiv_{\rm wtt}}
\newcommand{\bwttLE}{\le_{\bwtt}}
\newcommand{\bwttE}{\equiv_{\bwtt}}
\newcommand{\kwttLE}{\le_{\kwtt}}
\newcommand{\kwttE}{\equiv_{\kwtt}}
\newcommand{\nwttLE}{\le_{\nwtt}}
\newcommand{\nwttE}{\equiv_{\nwtt}}
\newcommand{\ttL}{<_{\rm tt}}
\newcommand{\kttLE}{\le_{\ktt}}
\newcommand{\kttL}{<_{\ktt}}
\newcommand{\kttE}{\equiv_{\ktt}}
\newcommand{\npttLE}{\le_{\ntt}}
\newcommand{\nttL}{<_{\ntt}}
\newcommand{\nttE}{\equiv_{\ntt}}
\newcommand{\hpTLE}{\le_{\rm T}^{\rm h}}
\newcommand{\NTLE}{\not\le_{\rm T}}
\newcommand{\TL}{<_{\rm T}}
\newcommand{\mL}{<_{\rm m}}
\newcommand{\pTL}{<_{\rm T}^{\rm p}}
\newcommand{\pT}{\le_{\rm T}^{\rm p}}
\newcommand{\pTLE}{\le_{\rm T}^{\rm p}}
\newcommand{\pttE}{\equiv_{{\rm tt}}^p}
\newcommand{\doub}{\es''}
\newcommand{\trip}{\es'''}
%
%
\newcommand{\inter}{\cap}
\newcommand{\union}{\cup}
\newcommand{\sig}[1]{\sigma_{#1} }
\newcommand{\s}[1]{\s_{#1}}
\newcommand{\LMA}{{\rm L}(M^A)}
\newcommand{\Ah}{{\hat A}}
\newcommand{\monus}{\;\raise.5ex\hbox{{${\buildrel
    \ldotp\over{\hbox to 6pt{\hrulefill}}}$}}\;}
\newcommand{\dash}{\hbox{-}}
\newcommand{\infinity}{\infty}
\newcommand{\ie}{\hbox{ i.e.  }}
\newcommand{\eg}{\hbox{ e.g.  }}

%
%
%
%
%
%
\newcounter{savenumi}
\newenvironment{savenumerate}{\begin{enumerate}
\setcounter{enumi}{\value{savenumi}}}{\end{enumerate}
\setcounter{savenumi}{\value{enumi}}}
\newtheorem{theoremfoo}{Theorem}[section] 
\newenvironment{theorem}{\pagebreak[1]\begin{theoremfoo}}{\end{theoremfoo}}
\newenvironment{repeatedtheorem}[1]{\vskip 6pt
\noindent
{\bf Theorem #1}\ \em
}{}

\newtheorem{lemmafoo}[theoremfoo]{Lemma}
\newenvironment{lemma}{\pagebreak[1]\begin{lemmafoo}}{\end{lemmafoo}}
\newtheorem{conjecturefoo}[theoremfoo]{Conjecture}
\newtheorem{research}[theoremfoo]{Line of Research}
\newenvironment{conjecture}{\pagebreak[1]\begin{conjecturefoo}}{\end{conjecturefoo}}

\newtheorem{conventionfoo}[theoremfoo]{Convention}
\newenvironment{convention}{\pagebreak[1]\begin{conventionfoo}\rm}{\end{conventionfoo}}

\newtheorem{porismfoo}[theoremfoo]{Porism}
\newenvironment{porism}{\pagebreak[1]\begin{porismfoo}\rm}{\end{porismfoo}}

\newtheorem{gamefoo}[theoremfoo]{Game}
\newenvironment{game}{\pagebreak[1]\begin{gamefoo}\rm}{\end{gamefoo}}

\newtheorem{corollaryfoo}[theoremfoo]{Corollary}
\newenvironment{corollary}{\pagebreak[1]\begin{corollaryfoo}}{\end{corollaryfoo}}

\newtheorem{claimfoo}[theoremfoo]{Claim}
\newenvironment{claim}{\pagebreak[1]\begin{claimfoo}}{\end{claimfoo}}

\newtheorem{openfoo}[theoremfoo]{Open Problem}
\newenvironment{open}{\pagebreak[1]\begin{openfoo}\rm}{\end{openfoo}}

\newtheorem{exercisefoo}{Exercise}
\newenvironment{exercise}{\pagebreak[1]\begin{exercisefoo}\rm}{\end{exercisefoo}}

\newcommand{\fig}[1] 
{
 \begin{figure}
 \begin{center}
 \input{#1}
 \end{center}
 \end{figure}
}

\newtheorem{potanafoo}[theoremfoo]{Potential Analogue}
\newenvironment{potana}{\pagebreak[1]\begin{potanafoo}\rm}{\end{potanafoo}}

\newtheorem{notefoo}[theoremfoo]{Note}
\newenvironment{note}{\pagebreak[1]\begin{notefoo}\rm}{\end{notefoo}}

\newtheorem{notabenefoo}[theoremfoo]{Nota Bene}
\newenvironment{notabene}{\pagebreak[1]\begin{notabenefoo}\rm}{\end{notabenefoo}}

\newtheorem{nttn}[theoremfoo]{Notation}
\newenvironment{notation}{\pagebreak[1]\begin{nttn}\rm}{\end{nttn}}

\newtheorem{empttn}[theoremfoo]{Empirical Note}
\newenvironment{emp}{\pagebreak[1]\begin{empttn}\rm}{\end{empttn}}

\newtheorem{examfoo}[theoremfoo]{Example}
\newenvironment{example}{\pagebreak[1]\begin{examfoo}\rm}{\end{examfoo}}

\newtheorem{dfntn}[theoremfoo]{Def}
\newenvironment{definition}{\pagebreak[1]\begin{dfntn}\rm}{\end{dfntn}}

\newtheorem{propositionfoo}[theoremfoo]{Proposition}
\newenvironment{proposition}{\pagebreak[1]\begin{propositionfoo}}{\end{propositionfoo}}
\newenvironment{prop}{\pagebreak[1]\begin{propositionfoo}}{\end{propositionfoo}}

\newenvironment{proof}
    {\pagebreak[1]{\narrower\noindent {\bf Proof:\quad\nopagebreak}}}{\QED}
\newenvironment{sketch}
    {\pagebreak[1]{\narrower\noindent {\bf Proof sketch:\quad\nopagebreak}}}{\QED}
\newenvironment{comment}{\penalty -50 $(*$\nolinebreak\ }{\nolinebreak $*)$\linebreak[1]\ }

\newenvironment{algorithm}[1]{\bigskip\noindent ALGORITHM~#1\renewcommand{\theenumii}{\arabic{enumii}}\renewcommand{\labelenumii}{Step \theenumii :}\begin{enumerate}}{\end{enumerate}END OF ALGORITHM\bigskip}

\newenvironment{protocol}[1]{\bigskip\noindent PROTOCOL~#1\renewcommand{\theenumii}{\arabic{enumii}}\renewcommand{\labelenumii}{Step \theenumii :}\begin{enumerate}}{\end{enumerate}END OF PROTOCOL\bigskip}

\newenvironment{red}[1]{\noindent REDUCTION~#1\renewcommand{\theenumii}{\arabic{enumii}}\renewcommand{\labelenumii}{Step \theenumii :}\begin{enumerate}}{\end{enumerate}END OF REDUCTION}

\newenvironment{con}{\noindent CONSTRUCTION\renewcommand{\theenumii}{\arabic{enumii}}\renewcommand{\labelenumii}{Step \theenumii :}\begin{enumerate}}{\end{enumerate}END OF CONSTRUCTION}

\newenvironment{alg}[1]{\bigskip\noindent ALGORITHM~#1\renewcommand{\theenumii}{\arabic{enumii}}\renewcommand{\labelenumii}{Step \theenumii :}\begin{enumerate}}{\end{enumerate}END OF ALGORITHM\bigskip}

\newcommand{\yyskip}{\penalty-50\vskip 5pt plus 3pt minus 2pt}
\newcommand{\blackslug}{\hbox{\hskip 1pt
        \vrule width 4pt height 8pt depth 1.5pt\hskip 1pt}}
\newcommand{\QED}{{\penalty10000\parindent 0pt\penalty10000
        \hskip 8 pt\nolinebreak\blackslug\hfill\lower 8.5pt\null}
        \par\yyskip\pagebreak[1]}

\newcommand{\BBB}{{\penalty10000\parindent 0pt\penalty10000
        \hskip 8 pt\nolinebreak\hbox{\ }\hfill\lower 8.5pt\null}
        \par\yyskip\pagebreak[1]}
     
\newcommand{\PYI}{CCR-8958528}
\newtheorem{factfoo}[theoremfoo]{Fact}
\newenvironment{fact}{\pagebreak[1]\begin{factfoo}}{\end{factfoo}}
\newenvironment{acknowledgments}{\par\vskip 6pt\footnotesize Acknowledgments.}{\par}




\newenvironment{construction}{\bigbreak\begin{block}}{\end{block}
    \bigbreak}

\newenvironment{block}{\begin{list}{\hbox{}}{\leftmargin 1em
    \itemindent -1em \topsep 0pt \itemsep 0pt \partopsep 0pt}}{\end{list}}


\dimen15=0.75em
\dimen16=0.75em

\newcommand{\lblocklabel}[1]{\rlap{#1}\hss}

\newenvironment{lblock}{\begin{list}{}{\advance\dimen15 by \dimen16
    \leftmargin \dimen15
    \itemindent -1em
    \topsep 0pt
    \labelwidth 0pt
    \labelsep \leftmargin
    \itemsep 0pt
    \let\makelabel\lblocklabel
    \partopsep 0pt}}{\end{list}}


\newenvironment{lconstruction}[2]{\dimen15=#1 \dimen16=#2
  \bigbreak\begin{block}}{\end{block}\bigbreak}

\newcommand{\Comment}[1]{{\sl ($\ast$\  #1\ $\ast$)}}


\newcommand{\Kobler}{K\"obler}
\newcommand{\Schoning}{Sch\"oning}
\newcommand{\Toran}{Tor\'an}
\newcommand{\Balcazar}{Balc{\'a}zar}
\newcommand{\Diaz}{D\'{\i}az}
\newcommand{\Gabarro}{Gabarr{\'o}}
\newcommand{\Laszlo}{L{\'a}szl{\'o}}
\newcommand{\Erdos}{Erd\"os }
\newcommand{\Erdosns}{Erd\"os}
\newcommand{\Sarkozy}{S{\'a}rk{\"o}zy }
\newcommand{\Sarkozyns}{S{\'a}rk{\"o}zy}
\newcommand{\Szekely}{Sz{\'e}kely }
\newcommand{\Szekelyns}{Sz{\'e}kely}
\newcommand{\Szemeredi}{Szemer{\'e}di }
\newcommand{\Szemeredins}{Szemer{\'e}di}
\newcommand{\Holder}{H{\"o}lder }
\newcommand{\Holderns}{H{\"o}lder}
\newcommand{\Chv}{Chv{\'a}tal }
\newcommand{\Chvns}{Chv{\'a}tal}

\newcommand{\CETppn}{\C {(p+1)^n} {\es'''}}
\newcommand{\logcp}{\ceil{\log(c+1)}}
\newcommand{\plogc}{\ceil{  {{\log(c)}\over{\log(p+1)}}      }   }
\newcommand{\plogcM}{\ceil{  {{\log(c)}\over{\log(p+1)}}      }-1   }
\newcommand{\plogcp}{\ceil{  {{\log(c+1)}\over{\log(p+1)}}      }   }
\newcommand{\plogcpm}{\ceil{  {{\log(c+1)}\over{\log(p+1)}}      }-1   }
\newcommand{\plogcPM}{\ceil{  {{\log(c)}\over{\log(p+1)}}      }-1   }
\newcommand{\plogbp}{\ceil{\log_{p+1}(b+1)}}
\newcommand{\plogbpm}{\ceil{\log_{p+1}(b+1)}-1}
\newcommand{\logc}{\ceil{\log(c)}}
\newcommand{\logcpm}{\ceil{\log(c+1)}-1}
\newcommand{\FQplogcpest}{\FQ{\plogcp}{\C p \trip}}
\newcommand{\FQplogcpmX}{\FQ{\plogcp -1}{\C p X}}
\newcommand{\FQplogcpmest}{\FQ{\plogcp-1}{\C p \trip}}
\newcommand{\FQplogcest}{\FQ{\plogc}{\C p \trip}}
\newcommand{\FQplogcMest}{\FQ{\plogcM}{\C p \trip}}
\newcommand{\FQplogcpm}{\FQ{\plogcp -1}{\C p K}}
\newcommand{\FQplogcp}{\FQ{\plogcp}{\C p K}}
\newcommand{\Gehrs}{G_{s}}
\newcommand{\Gehr}{G}
\newcommand{\Gers}{G_{e,s}^r}
\newcommand{\Ger}{G_e^r}
\newcommand{\Ge}{G_e}
\newcommand{\Ghrfk}{G_{f_{k,k+1}(x)}^{hr}}
\newcommand{\Ghrf}{G_{f_{k,m}(x)}^{hr}}
\newcommand{\Gn}{G^n}
\renewcommand{\L}[1]{{\hat L}_{#1}}
\newcommand{\Lei}{L_{\ei}}
\newcommand{\Lsei}{L_{\ei}^s}
\newcommand{\Ve}{V_e}
\newcommand{\Xei}{X_{\ang{e,i}}}
\newcommand{\cGehr}{\chi(G)}
\newcommand{\cGehrs}{\chi(G_{s})}
\newcommand{\ceiling}[1]{\lc{#1}\rc}
\newcommand{\chir}{\chi^r}
\newcommand{\eoet}{\ang {e_1,e_2} }
\newcommand{\logbm}{\ceiling{\log (b-1)}}
\newcommand{\logbpm}{\logbp -1}
\newcommand{\logbp}{\ceiling{\log (b+1)}}
\newcommand{\logb}{\ceiling{\log b}}
\newcommand{\logcGehr}{\log \cGehr }
\newcommand{\logcm}{\lc \log c \rc-1}
\newcommand{\logcppm}{\lc {\log (c+1)\over p} \rc-1}
\newcommand{\logcpp}{\lc {\log (c+1)\over p} \rc}
\newcommand{\lognp}{\lc \log (n+1) \rc}
\newcommand{\logpp}{\log(p+1)}
\newcommand{\logrcGehr}{\log \rcGehr }
\newcommand{\meis}{m_{\ei}^s}
\newcommand{\mei}{m_{\ang{e,i}}}
\newcommand{\plogbm}{\ceiling{\log (b-1)\over\logpp}}
\newcommand{\plogb}{\lc {\log b\over \log (p+1)} \rc}
\newcommand{\plogcmM}{\ceiling{\log(c-1)\over\log(p+1)}-1}
\newcommand{\plogcm}{\ceiling{\log (c-1)\over\logpp}}
\newcommand{\plogcpM}{\lc{\log(c+1)\over\logpp}\rc-3}
\newcommand{\pplogcmM}{\lc{\log(c-1)-1\over p}\rc}
\newcommand{\pplogcM}{\lc{\log(c)-1\over p}\rc}
\newcommand{\pplogcm}{\lc{\log(c-1)\over p}\rc}
\newcommand{\pplogcpm}{\lc{\log(c+1)-1\over p}\rc}
\newcommand{\pplogcp}{\lc{\log(c+1)\over p}\rc}
\newcommand{\pplogcpM}{\lc{\log(c+1)\over p}\rc-1}
\newcommand{\qam}{q^{a-1}}
\newcommand{\rcGehr}{\chi^r(G)}
\newcommand{\rchi}{\chi^r}
\newcommand{\ttlogcm}{2^{\logcm}}
\newcommand{\ttlogcpm}{2^{\logcpm}}
\newcommand{\unioneinf}{\bigcup_{e=0}^{\infty}}
\newcommand{\unionsinf}{\bigcup_{s=0}^{\infty}}
\newcommand{\logsum}{{\rm logsum}}

\newcommand{\fcG}{f(\chi(G))}
\newcommand{\gcG}{g(\chi(G))}
\newcommand{\fcrG}{f(\chir(G))}
\newcommand{\gcrG}{g(\chir(G))}
\newcommand{\fpcG}{f_p(\chi(G)}
\newcommand{\gpcG}{g_p(\chi(G)}
\newcommand{\pfcG}{\fcG\over p}
\newcommand{\pgcG}{\gcG\over p}
\newcommand{\pfcrG}{\fcrG\over p}
\newcommand{\pgcrG}{\gcrG\over p}
\newcommand{\gibar}{\overline{GI}}
\newcommand{\adv}{{\rm ADV}}

\centerline{\bf Low,  Superlow, and Superduperlow Sets:}
\centerline{\bf An Exposition of a Known But Not Well-Known Result}
\centerline{By William Gasarch}

\section{Introduction}

\begin{notation}~
\begin{enumerate}
\item
$M_0, M_1,\ldots$ is be a standard list of Turing Machines.
\item
$W_e$ is the domain of $M_e$. Hence
$W_0, W_1, \ldots$ is a list of all c.e.\ sets.
\item
$M_0^{()}, M_1^{()},\ldots$ is a standard list of oracle Turing Machines.
\item
$K$ is the set $\{e \st M_e(e)\cvg \}$.
\item
If $A$ is a set then $A'= \{e \st M_e^A(e)\cvg \}$.
\item
A set $A$ is {\bf Low} if $A'\TLE K$.
\item
A set $A$ is {\bf Superlow} if $A'\le_{tt} K$.
\item
A set $A$ is {\bf Superduperlow} if $A'\le_{btt} K$.
\end{enumerate}
\end{notation}

By a finite injury priority argument one can construct a noncomputable low c.e.\ set.
On closer examination of the proof you can extract that $A$ is actually superlow.
(We include this proof in the appendix.)
This raises the question: Is there a noncomputable superduperlow set $A$?

I asked about this at a logic conference and found out:

\begin{enumerate}
\item
Four prominent computability theorists thought that there was no such set; however,
none knew of a proof or reference.
\item
Carl Jockusch, also a prominent computability theorist, knew of three 
unpublished proofs
(one by Bickford and Mills, one by Phillips, and one by himself) 
and also a more complicated published proof by Mohrerr~\cite{superduperlow}.
She actually proved the stronger result that if 
$A^{tt} \leq_{btt} B'$ then $A \TLE B$,
as did Bickford and Mills.
\end{enumerate}

In this manuscript we (1) give the unpublished proof that is due to Jockusch,
and (2) give a new proof by Frank Stephan.

We will use the following Lemma, called the Shoenfield Limit Lemma.

\begin{lemma}\label{le:sh}
$A\TLE K$ iff there exists a computable function $h:\nat\times\nat \into \nat$ such that
$$A(x) = \lim_{s\goes\infinity} h(x,s).$$
\end{lemma}

\section{Bi-immune and Hyperimmune Sets}

\begin{definition}~
\begin{enumerate}
\item
A set $C$ is {\it immune} if $C$ is infinite and has no infinite c.e.\ subsets
\item
A set $C$ is {\it bi-immune} if both $C$ and $\overline{C}$ are immune.
\item
If $B$ is a set then {\it the principal function of $B$}, denoted $p_B$, is
defined as 
$$p_B(x)=\hbox{ the $x^{\rm th}$ element of $B$}.$$
\item
A function $f$ is {\it majorized} by a function $g$ if,
for all $x$, $f(x)<g(x)$.
\item
A set $B$ is {\it hyperimmune} if  $p_B$ 
is not majorized by any computable function.
(See note for why this is called {\it hyperimmune}.)
\end{enumerate}
\end{definition}

\begin{note}
There is a different (but equivalent) definition of hyperimmune that
is an extension of immune (hence the name).
We present it here even though we will not use it.
Recall that 
$D_n$ is the set corresponding to the bit vector of $n$ in binary.
A set $A$ is {\it hyperimmune} if $\overline{A}$ is infinite
and there is no computable 
function $i$ such that the following occurs:
\begin{itemize}
\item
The sets $D_{i(1)}, D_{i(2)}, \ldots$ are disjoint.
\item
For all $m$ $D_{i(m)}\cap A \ne\es$.
\end{itemize}
\end{note}

The following Lemma is due to Miller and Martin~\cite{MillerMartin}.
We include the proof since the original article is behind a paywall and
hence lost to humanity forever.

\begin{lemma}\label{le:existshyper}
If $\es <_T A\TLE K$ then
there exists a hyperimmune set $B$ such that $B\TLE A$.
\end{lemma}

\begin{proof}

Since $A\TLE K$ there exists, by Lemma~\ref{le:sh},
a computable $h$ such that

$$A(x) = \lim_{s\goes\infinity} h(x,s).$$

Let

$$
f(x) = \hbox{ the least $s\ge \max\{x,f(x-1)\}$ such that }
 (\forall y\le x)[A(y)=h(y,s)].
$$

Let $B$ be the image of $f$.  
Clearly $B\TLE A$.
Note that $f=p_B$, the principal function of $B$.
Hence it suffices to show that
$f$ cannot be majorized by any computable function.

Assume, by way of contradiction, that there is a computable $g$ such that,
for all $x$, $f(x)<g(x)$.  We use this to obtain an algorithm for $A$.
Given $x$, we want to determine $A(x)$.

Note that, for all $y$, 
$$y \le f(y) \le g(y)$$

Let $y\ge x$. Imagine what would happen if 

$$h(x,y)=h(x,y+1)=\cdots=h(x,f(y))=\cdots=h(x,g(y)) = b.$$

$$A(x) = h(x,f(y))= b.$$

Hence we would know $A(x)$. Therefore we need to find such a $y$.
If we knew that one existed we could just look for it.

One does exist! Let $y$ be such that 

$$h(x,y)=h(x,y+1)=\cdots=$$

Such a $y$ exists since $h$ reaches a limit.
This $y$ clearly suffices.
We cannot find this particular $y$ but we do not need to.
We need only find a $y$ such that
$$h(x,y)=h(x,y+1)=\cdots=h(x,g(y))  = b.$$

Here is the formal algorithm for $A$.

\begin{enumerate}
\item
Input($x$)
\item
Find a $y\ge x$ such that 
$$h(x,y)=h(x,y+1)=\cdots=h(x,g(y))  = b$$
\item
Output the value $b$.
\end{enumerate}

Thus $A$ is computable--- a contradiction.
Hence $f$ cannot be majorized by any computable function.

Therefore we have a set $B\TLE A$ such that $B$ is hyperimmune.

\end{proof}

The following is a result of Carl Jockusch~\cite{jockbi}

\begin{lemma}\label{le:existsbi}
For every hyperimmune $B$ there exists a bi-immune $C\TLE B$.
\end{lemma}

\begin{proof}
Let  $B$  be hyperimmune. Let  $f=p_B$, the principal function of $B$.
Since $B$ is hyperimmune $f$ is not 
majorized by any computable function.
We use this to construct a bi-immune $C\leq_T B$.
To ensure that $C$ is bi-immune we make sure that 
$C$ satisfies the following requirements:

$$R_e: W_e \hbox{ infinite } \rightarrow (W_e\cap C\ne \es \wedge W_e\cap\overline{C}\ne\es).$$

\noindent
CONSTRUCTION

\noindent
{\it Stage 0:} For all $e$, $R_e$ is not satisfied.

\medskip

\noindent
{\it Stage s:} Find the least $e\le s$, if it exists, such that
$R_e$ is not satisfied and $W_{e, f(s)}$ has at least two elements $x_1,x_2\ge s$ which have not
yet been put into  $C$  or $\overline{C}$.  
Put $x_1$ into $C$, $x_2$ into $\overline{C}$,
and declare  $R_e$ {\it satisfied}.
(it will never become unsatisfied).
We also say that $R_e$ has {\it acted}.
If there is no such  $e$,  do nothing.

\noindent
END OF CONSTRUCTION

We have  $C \leq_T  B$  since the construction is  $B$-computable
and  $C(n)$  is decided by stage  $n$  .   (For definiteness, a number is
in  $C$  iff the construction puts it into  $C$ .)

We show that $C$ is bi-immune by showing that it satisfies each requirement.
We assume that $R_1,\ldots,R_{e-1}$ are satisfied and show that $R_e$ is satisfied.
There are two cases.
\begin{enumerate}
\item
$W_e$ is finite.  Then clearly $R_e$ is satisfied. 
\item
$W_e$ is infinite.
Assume, by way of contradiction, that $W_e$ is not satisfied.
From this we will construct a computable function $g$ that majorizes $f$ which
will be the contradiction.
Let $s_0$ be such that by state $s_0$ all of $R_1,\ldots,R_{e-1}$ that are going to act have
acted. So for all $s\ge s_0$ $R_e$ is not satified yet fails to act! Why!?

Let  $g(s)$  be the least  $t\le s_0$  such that $W_{e,t}$
has at least $2s+2$  elements  $\geq  s$.   If $g(s)\le f(s)$ then in stage $s$,
since at most $2s$ elements have been determined, there must exist $x_1,x_2\in W_{e,t}$
And yet $R_e$ has not acted! Why not? It must be that 
$(\forall s\ge s_0)[g(s) > f(s)]$.
And $g$ is computable! Hence there is a computable function that majorizes
$f$. This is a contradiction.
\end{enumerate}
\end{proof}

\begin{lemma}\label{le:all}
Let  $\es \TLE A \TLE K$.
\begin{enumerate}
\item
If $A$ is not computable then there exists $C$ bi-immune such that $C\TLE A$.
\item
If there is no bi-cimmune set $C\TLE A$ then $A$ is computable (this is just the
contrapositive of part 1 so we won't proof it.)
\end{enumerate}
\end{lemma}

\begin{proof}
1) By Lemma~\ref{le:existshyper} there is a hyperimmune set $B\TLE A$.
By Lemma~\ref{le:existsbi} there is a bi-immune set $C\TLE B$. Hence
there is a bi-immune set $C\TLE B \TLE K$.
\end{proof}

\section{Final Theorem}

\begin{definition}
A set $D$ is {\it weakly $n$-c.e.} if there
exists a function $h$ such that
\begin{itemize}
\item
$D(x) = \lim_{s\goes\infinity} h(x,s)$
\item
$|\{s \st h(x,s)\ne h(x,s+1)\}|\le n$.
\end{itemize}
\end{definition}

The following easy lemma we leave to the reader.

\begin{lemma}\label{le:nce}
If $D\le_{btt} K$ then there exists an $n$ such that
$D$ is $n$-c.e.\ but not $(n-1)$-c.e.
\end{lemma}

The following is an unpublished proof of Jockush.

\begin{theorem}
If $A$ is superduperlow then $A$ is decidable.
\end{theorem}

\begin{proof}
Since $A$ is superduperlow $A'\le_{btt} K$.

Let $D\TLE A$.  We will show that $D$ is not bi-immune.
Note that 

$$D\le_m A' \le_{btt} K.$$

Hence $D\le_{btt} K$.
By Lemma~\ref{le:nce} there is some $n$ such that $D$ is weakly $n$-c.e.\ 
but not weakly $(n-1)$-c.e.
Let $h$ be such that

\begin{itemize}
\item
$D(x) = \lim_{s\goes\infinity} h(x,s)$
\item
$|\{s \st h(x,s)\ne h(x,s+1)|\le n$.
\end{itemize}

Let 

$$E= \{ x \st |\{s \st h(x,s)\ne h(x,s+1)|=n\}\}.$$

$E$ is infinite, else $D$ is weakly $(n-1)$-c.e.
Let 

$$E_0= \{ x \st h(x,0)=0 \wedge |\{s \st h(x,s)\ne h(x,s+1)|=n\}\}.$$

$$E_1= \{ x \st h(x,0)=1 \wedge |\{s \st h(x,s)\ne h(x,s+1)|=n\}\}.$$

$E_0$ and $E_1$ are both c.e:

\bigskip

$E_0=\{ x \st h(x,0)=0\wedge (\exists s_1<\cdots<s_{n+1})(\forall i\le n)[h(x,s_i)\ne h(x,s_{i+1})] \}.$

\bigskip

$E_1=\{ x \st h(x,0)=1\wedge (\exists s_1<\cdots<s_{n+1})(\forall i\le n)[h(x,s_i)\ne h(x,s_{i+1})] \}.$

At least one of $E_0$ or $E_1$ is infinite. We assume it is $E_0$
(the case for $E_1$ is similar). We also assume that $n$ is even
(the case of $n$ odd is similar).
$E_0$ is an infinite c.e.\ subset of $D$.
In all of the omitted cases you either get $D$ has an infinite c.e.\ subset or $\overline{D}$ has an infinite c.e.\ subset.
Hence $D$ is not bi-immune.

The upshot is that for {\it every} set $D\TLE A$ is not bi-immune.
By Lemma~\ref{le:all} $A$ is computable.
\end{proof}

\section{Another Proof}

In this section we present a proof by Frank Stephan
that uses concepts from Bounded Queries.

The following theorem was first proven by Kummer~\cite{kummer}.
See also the survey of bounded queries in computability theory
by Gasarch~\cite{gems} for a different proof which is free and online.

\begin{definition}
Let $A\subseteq \nat$ and $n\in \nat$.
\begin{enumerate}
\item
$\mem n A : \nat^n \into \bits n$ is the following function:
$$\mem n A \xtu n = A(x_1)\cdots A(x_n).$$
\item
$\card n A : \nat^n \into \nat$ is the following function:
$$\card n A \xtu n = |\xset n \cap A|.$$
\item
A function $f$ is in $\EN m$ if there exists $m$ partial computable
functions $f_1,\ldots,f_m$ such that
$$(\forall x)(\exists i)[f(x)=f_i(x)].$$
We will use the following equivalent definition:
there is a process that will, given $x$, enumerate at most $m$ ouptuts
one of which is the answer.
\end{enumerate}
\end{definition}

Clearly, for all $A$,  $\card n A \in \EN {n+1}$.
Kummer's theorem states that, for undecidable sets,  this is the best one can do.

\begin{theorem}\label{th:card}
For all $A$, if $\card n A \in \EN n$ then $A$ is computable.
\end{theorem}

We use Theorem~\ref{th:card} to show that all superduperlow sets are
decidable.

We need an easy lemma

\begin{lemma}\label{le:memK}
$$\mem n K \in \EN {n+1}.$$
\end{lemma}

\begin{theorem}
If $A$ is superduperlow then $A$ is decidable.
\end{theorem}

\begin{proof}
Assume that $A$ is superduperlow.
Let $k$ be such that $A \ktt K$.
Note that $k$ is a constant.

We show that for some (large enough) $n$, $\card {2^n-1} A \in \EN {2^n-1}$.
We will choose $n$ later.

Let $A_1,\ldots,A_n$ be the following sets.

$$A_i = \{ \xtu {2^n-1} \st \hbox{ the $i$th bits of $\card {2^n-1} A$ is 1 } \}.$$

For each $i$ $A_i \le_{T} A$. Hence $A_i \le_m A' \le_{k-tt} K$.
Therefore $\card {2^n-1} A$ can be computed with $kn$ queries to $K$.
With this in mind we present the following procedure for
$$\card {2^n-1} A \in \EN {kn+1}.$$

\begin{enumerate}
\item
Input $\xtu {2^n-1}$.
\item
Using $A_i \le_{k-tt} K$ find $kn$ numbers $y_1,\ldots,y_{kn}$
such that if we knew $\mem {kn} A \ytu {kn} $ then we would know
$\xtu {2^n-1}$.
\item
By Lemma~\ref{le:memK} $\mem  {kn} K \in \EN {kn+1}$.
Run this enumeation algorithm.
Every time a candidate for 
$\mem {kn} K$ is enumerated, use it to obtain
a candidate for 
$\card {2^n-1} A.$
\end{enumerate}

By the above enumeation algorithm $\card {2^n-1} A \in \EN {kn+1}.$
Take $n$ large enough so that $kn+1 \le 2^n-1$ to obtain that
$\card A {2^n-1} \in \EN {2^n-1}$ and hence $A$ is computable.
\end{proof}

\appendix
\section{There exists an undecidable  c.e.\ Superlow Set}

\begin{theorem}
There exists an undecidable c.e.\ superlow set.
\end{theorem}

\begin{proof}
We construct a c.e. set $A$ that satisfies the following requirements:

$P_e: W_e\hbox{ infinite } \implies W_e \cap A \ne \es.$

These are called {\it positive requirements} since they act by putting numbers
into $A$. It is easy to show that if $A$ is co-infinite and all of the $P_e$'s are satisfied then $A$ is undecidable.
(We will also make $A$ co-infinite though we do not state it as a formal requirement.)

$N_e: (\IO s)[M_{e,s}^{A_s}(e)\cvg] \implies M_e^A(e)\cvg.$

These are called {\it negative requirements} since they will act by restraining
numbers from coming into $A$ in order to protect a computation from being
injured. Associated to every $N_e$ will be a restraint function $r(e,s)$.
This is $N_e$ saying {\it you cannot put an element into $A$ that is $\le r(e,s)$}.
This restraint will be respected by the lower priority positive requirements
($P_e$, $P_{e+1}$, etc.)
but not by the higher priority positive requirements
($P_0$, $P_1$, \ldots, $P_{e-1}$).

The requirements are in the following priority ordering

$$N_0, P_0, N_1, P_1, \ldots$$

\noindent
{\bf CONSTRUCTION}

\noindent
{\bf Stage 0:} $A_0=\es$, $(\forall e)[r(e,0)=0]$, for all $e$ $P_e$ is not satisfied.

\smallskip

\noindent
{\bf Stage $s$:} Visit each requirement in turn via the priority ordering, up to $P_s$.

\noindent
{\bf Case 1:} A positive requirement $P_e$. If 
(a) $P_e$ is not satisfied,
(b) there exists $x\in W_{e,s}$ such that $x\ge 2e$ and $x\ge \max_{i\le e} r(e,s)$ then $P_e$ acts
by putting $x$ into $A$. $P_e$ is declared satisfied. 
For every $i<e$ set $r(i,s)=0$.
(This is not really needed but makes the proof cleaner.)
We say that $N_i$ is {\it injured}.
Note that $P_e$ will never become unsatisfied.

\noindent
{\bf Case 2:} A negative requirement $N_e$. If 
$M_{e,s}^{A_s}(e)\cvg$ then set $r(e,s)$ to be the largest number that is queried
in this computation.

\noindent
{\bf END OF CONSTRUCTION}

\bigskip

\noindent
{\bf Claim 1:} Every $P_e$ acts finitely often.

\noindent
{\bf Proof of Claim 1:} Once $P_e$ acts it is satisfied and never acts again.

\bigskip

\noindent
{\bf End of Proof of Claim 1:}

\bigskip

\noindent
{\bf Claim 2:} For all $e$,  $\lim_{s\goes\infinity} r(e,s)$ is finite and
$N_e$ is satisfied.

\noindent
{\bf Proof of Claim 2:} Let $s_0$ be such that, for all $i<e$, $P_e$ will never
act past stage $s_0$. Note that $s_0$ exists by Claim 1.
Past stage $s_0$ $N_e$ will never get injured.
Hence if there exists $s>s_0$ such that $r(e,s)$ is set to a non-zero value
then it will remain there.
Note that $N_e$ will only be injured finitely often. Hence it is satisfied.

\noindent
{\bf End of Proof of Claim 2:}

\bigskip

\noindent
{\bf Claim 3:} Every $P_e$ is satisfied.

\noindent
{\bf Proof of Claim 3:} 
If $W_e$ is finite then $P_e$ is satisfied.
Hence we assume that $W_e$ is infinite.
Let $s_0$ be such that for all $i<e$ $\lim_{s\goes \infinity} r(i,s) = r(i,s_0)$.
Let $R(e)=\max_{i<e} r(i,s_0)$. Since $W_e$ is infinite there will
be an $x\ge \max\{2e,R(e)\}$ that is enumerated into $W_e$ at some stage $s>\max\{s_0,e\}$.
If $P_e$ is not yet satisfied then $P_e$ will act at stage $s$ and be satisfied.

\noindent
{\bf End of Proof of Claim 3:}

\bigskip

\noindent
{\bf Claim 4:} $A$ is co-infinite.

\noindent
{\bf Proof of Claim 4:} 
Look at the numbers $S_e=\{1,2,\ldots,2e,2e+1\}$. Since $P_{e+1}$ only uses numbers $\ge 2e+2$,
the only positive requirements that will use elements of $S_e$ are $P_0,\ldots,P_e$.
Hence at most $e+1$ of the elements of $S_e$ will enter $A$. Hence at least $e$ of the
elements of $S_e$ will not enter $A$. Since this is true for all $e$, $A$ is co-infinite.

\noindent
{\bf End of Proof of Claim 4:}

We now give the proof that $A$ is low, which is the standard conclusion of
this argument. We will then show that $A$ is actually superlow.
\noindent

\bigskip

\noindent
{\bf Claim 5:} $A$ is low.

\noindent
{\bf Proof of Claim 5:} 

Algorithm to determine if $e\in A'$ given access to $K$.

\begin{enumerate}
\item
Set $s_{YES}=s_{NO}=0$.
\item
Ask $K$ ``$(\exists s\ge s_{NO})[M_{e,s}^{A_s}\cvg]$?''
If the answer is NO then $e\notin A'$.
If the answer is YES then find the least such $s$  and call it $s_{YES}$.
\item
Ask $K$ ``$(\exists s>s_{YES})[M_{e,s}^{A_s}\dvg]$?''
If the answer is NO then $e\in A'$.
If the answer is YES then find the least such $s$ and call it $s_{NO}$.
Goto step 2.
\end{enumerate}

Since $N_e$ is satisfied this process must terminate.

\noindent
{\bf End of Proof of Claim 5}

\bigskip

\noindent
{\bf Claim 6:} $A$ is superlow.

\noindent
{\bf Proof of Claim 6:} 

Note that the only requirements that can injure $N_e$
are $P_0,P_1,\ldots,P_{e-1}$.
These requirements act at most once. Hence
$N_e$ is injured at most $e$ times.
We can determine $e\in A'$ by asking the following questions:
For each $i\le e$ ask the two questions to $K$.
\begin{itemize}
\item
Is $N_e$ injured at least $i$ times.
\item
Is there a stage $s$ that occurs after $N_e$ is injured $i$ times where
$M_{e,s}^{A_s}(e)\cvg$.
\end{itemize}

(Note that we are just asking questions- we are not actually running any machines
that might not halt.)

We leave it as an exercise to show that the answers to these questions suffice
to determine if $e\in A'$.
\end{proof}

\section{Acknowledgment}

I thank Carl Jockusch and Frank Stephan for supplying most of the material
for this manuscript.


\begin{thebibliography}{1}

\bibitem{gems}
W.~Gasarch.
\newblock Gems in the field of bounded queries.
\newblock In Cooper and Goncharov, editors, {\em Computability and Models},
  2003.
\newblock \url{http://www.cs.umd.edu/~gasarch/papers/papers.html}.

\bibitem{jockbi}
C.~Jockusch.
\newblock The degrees of bi-immune sets.
\newblock {\em Zeitschrift f{\"u}r logik and Grundlagen d. Math (Has changed
  its name to Math Logic Quarterly)}, 15:135--140, 1986.
\newblock The article is behind a paywall and hence lost to humanity forever.

\bibitem{kummer}
M.~Kummer.
\newblock A proof of {B}eigel's cardinality conjecture.
\newblock {\em Journal of Symbolic Logic}, 57(2):677--681, June 1992.
\newblock
  \url{http://www.jstor.org/action/showPublication?journalCode=jsymboliclogic}.

\bibitem{MillerMartin}
W.~Miller and D.~A. Martin.
\newblock The degree of hyperimmune sets.
\newblock {\em Zeitschrift f{\"u}r logik and Grundlagen d. Math (Has changed
  its name to Math Logic Quarterly)}, 14:159--166, 1968.
\newblock The article is behind a paywall and hence lost to humanity forever.

\bibitem{superduperlow}
J.~Mohrherr.
\newblock A refinement of ${\rm low}_n$ and ${\rm high}_n$ for the r.e.
  degrees.
\newblock {\em Zeitschrift f{\"u}r logik and Grundlagen d. Math (Has changed
  its name to Math Logic Quarterly)}, 32:5--12, 1986.
\newblock The article is behind a paywall and hence lost to humanity forever.

\end{thebibliography}

\end{document}